\documentclass[10pt]{article}

\usepackage{mathptmx}
\usepackage{eucal}
\usepackage{amsmath}
\usepackage{amscd}
\usepackage{amssymb}
\usepackage{amsthm}
\usepackage{xspace}
\usepackage[all,tips]{xy}
\usepackage[dvips]{graphicx}
\usepackage{verbatim}
\usepackage{syntonly}
\usepackage{hyperref}

\theoremstyle{plain}
\newtheorem{thm}{Theorem}[section]
\newtheorem{cor}[thm]{Corollary}
\newtheorem{prop}[thm]{Proposition}
\newtheorem{lemma}[thm]{Lemma}
\newtheorem{ques}{Question}

\theoremstyle{definition}

\newenvironment{pf}
{\begin{proof}} {\end{proof}}

\DeclareMathOperator{\Aut}{Aut} \DeclareMathOperator{\Out}{Out}
 \DeclareMathOperator{\Cl}{Cl}
\DeclareMathOperator{\SL}{SL} 
\DeclareMathOperator{\GL}{GL}

\DeclareMathOperator{\Nil}{Nil} \DeclareMathOperator{\Ad}{Ad}

 \DeclareMathOperator{\Semi}{Semi}

\newcommand{\vp}{\varphi}

\newcommand{\bdef}{\overset{\text{def}}{=}}

\newcommand{\al}{\alpha}
\newcommand{\be}{\beta}
\newcommand{\ga}{\gamma}
\newcommand{\Ga}{\Gamma}
\newcommand{\te}{\theta}
\newcommand{\la}{\lambda}
\newcommand{\La}{\Lambda}

\newcommand{\ol}{\overline}

\newcommand{\wh}{\widehat}

\newcommand{\iny}{\infty}
\newcommand{\tri}{\ensuremath{\triangle}}
\newcommand{\es}{\emptyset}
\newcommand{\co}{\ensuremath{\colon}}

\newcommand{\innp}[1]{\left< #1 \right>}
\newcommand{\abs}[1]{\left\vert#1\right\vert}
\newcommand{\set}[1]{\left\{#1\right\}}

\newcommand{\pr}[1]{\left( #1 \right) }
\newcommand{\su}{\subset}

\newcommand{\ba}{\bigcap}
\newcommand{\bop}{\bigoplus}

\newcommand{\lra}{\longrightarrow}
\newcommand{\op}{\oplus}

\newcommand{\B}[1]{\ensuremath{\mathbf{#1}}}

\newcommand{\Fr}[1]{\ensuremath{\mathfrak{#1}}}

\newcommand{\Q}{\ensuremath{\B{Q}}}
\newcommand{\R}{\ensuremath{\B{R}}}
\newcommand{\Z}{\ensuremath{\B{Z}}}
\newcommand{\C}{\ensuremath{\B{C}}}

\usepackage{fancyhdr}

\fancypagestyle{plain}

\begin{document}


\title{\textbf{Controlling manifold covers of orbifolds}}
\author{D. B. McReynolds\thanks{Partially supported by an NSF postdoctoral fellowship.}}
\maketitle


\begin{abstract}
\noindent In this article we prove a generalization of
Selberg's lemma on the existence of torsion free, finite index
subgroups of arithmetic groups. Some of
the geometric applications are the resolution a conjecture of
Nimershiem and answers to questions of Long--Reid and the author.
\end{abstract}

\section{Introduction}\label{Intro}

\noindent For a compact orbifold $M$, there is no reason to expect
$M$ to possess a finite manifold cover. Indeed, even the existence of
a finite orbifold cover cannot be guaranteed. However, when
$\pi_1^{orb}(M)$ admits a faithful linear representation,
Selberg's lemma (see for instance \cite{Borel}) furnishes $M$ with many finite
manifold covers. Given their prolificacy, one might ask more geometrically of
these covers. Explicitly, we ask the following pair of questions.

\begin{ques}
\textsl{If $N$ is a properly immersed, $\pi_1$--injective submanifold of $M$,
can $N$ be lifted to an embedded submanifold in a finite manifold
cover of $M$?}
\end{ques}

\begin{ques}
\textsl{If $N$ is an immersed, totally geodesic submanifold of $M$, can $N$ to
be lifted to an embedded submanifold in a finite manifold cover of
$M$?}
\end{ques}

\noindent This article aims at resolving the first question in some
special cases---in the final section, we partially address the second question. Even in these special situations, there are some new
geometric applications. We have elected to postpone
the motivation for these geometric results until
Section \ref{GeometricAppsSect}. Before describing them, we give an abbreviated account of the associated algebraic problem. \smallskip\smallskip

\noindent The enterprize of promoting immersions to embeddings in
finite covers has received some attention in recent years. The
associated algebraic problem for the subgroup $\pi_1(N)$ of
$\pi_1^{orb}(M)$ is directly related to subgroup separability
(see \cite{Scott}). In this vein, we proved in \cite{McR1} a
result that promotes  $\pi_1$--injective immersions
to embeddings in finite covers when $M$ is arithmetic and $N$ is
infranil. Our present goal is to ensure the cover constructed in \cite{McR1}
can be taken to be a manifold. Algebraically, this requires a torsion free,
finite index subgroup $\La_0$ of $\pi_1^{orb}(M)$ that
contains $\pi_1(N)$. The main result of this article is the
resolution of this problem---throughout the remainder of this article, $[\eta]$ will denote the
$\GL(n,\Z)$--conjugacy class of an element $\eta$ of $\GL(n,\Z)$.

\begin{thm}\label{SemisimpleAvoid}
Let $\eta \in \GL(n,\Z)$ be a semisimple element and $\Ga<\GL(n,\Z)$
a torsion free, virtually unipotent subgroup. Then there exists a
finite index subgroup $\La_0$ of $\GL(n,\Z)$ such that $\Ga<\La_0$ and
$[\eta]\cap \La_0 = \es$.
\end{thm}

\begin{cor}\label{UnipotentSelbergLemma}
Let $\Gamma<\GL(n,\Z)$ be a torsion free, virtually unipotent subgroup.
Then there exists a torsion free, finite index subgroup
$\La_0$ of $\GL(n,\Z)$ such that $\Gamma < \La_0$.
\end{cor}

\noindent The main geometric application of Corollary
\ref{UnipotentSelbergLemma} given here is on the structure of cusp
cross sections of arithmetic orbifolds and manifolds. Specifically,
using the aforementioned subgroup separability result \cite[Theorem
3.1]{McR1} in tandem with Corollary \ref{GenUnipotentSelbergLemma}
(see Section \ref{GenSect}), we can promote $\pi_1$--injective
immersions of infranil manifolds into arithmetic orbifolds to
embeddings in some finite manifold cover of the target orbifold.
With this and our previous work in \cite{McR1, McR2}, we can derive
a few new geometric results. The first verifies a conjecture of
Nimershiem \cite[Conjecture 2']{Nim1}.

\begin{thm}[Nimershiem's conjecture]\label{Nimershiem}
Let $M$ be a closed flat $n$--manifold. Then the space of
similarity classes of flat structures that can be realized in cusp
cross sections of (arithmetic) hyperbolic $(n+1)$--manifolds is dense in the
space of flat similarity classes.
\end{thm}

\noindent This was previously known only for $n=1,2$, and $3$ (see
\cite{Nim1}). The following corollary of Theorem \ref{Nimershiem}
was also previously unknown.

\begin{cor}\label{Arith}
Every closed flat $n$--manifold is diffeomorphic to a cusp cross
section of an arithmetic hyperbolic $(n+1)$--manifold.
\end{cor}

\noindent Corollary \ref{Arith} upgrades the main result of
Long--Reid \cite[Theorem 1.1]{LR} to manifolds, answering a question
implicitly asked by Long and Reid \cite[p. 286]{LR}
(Nimershiem \cite[Conjecture 1']{Nim1} also conjectured this without an arithmetic
assumption). \smallskip\smallskip

\noindent Our next result is the extension of Theorem
\ref{Nimershiem} to the complex and quaternionic hyperbolic settings
via \cite[Theorem 3.5]{McR2} and Corollary
\ref{GenUnipotentSelbergLemma}.

\begin{thm}\label{Complex}

\begin{itemize}
\item[(a)]

Let $N$ be a closed almost flat manifold modeled on the Heisenberg
group $\Fr{N}_{2n+1}$. Then the space of similarity classes of
almost flat metrics on $N$ that can be realized in cusp cross
sections of complex hyperbolic $(n+1)$--manifolds is either empty or
dense in the space of almost flat metrics.

\item[(b)]

Let $N$ be a closed almost flat manifold modeled on the
quaternionic Heisenberg group $\Fr{N}_{4n+3}$. Then the space of
similarity classes of almost flat metrics on $N$ that can be
realized in cusp cross sections of quaternionic hyperbolic
$(n+1)$--manifolds is either empty or dense in the space of almost
flat metrics.
\end{itemize}
\end{thm}

\noindent In \cite[Theorem 5.4]{McR1}, we gave a necessary and sufficient condition on when
this set is non-empty. This provides the following corollary which answers a question asked
in \cite[Section 8]{McR1}.

\begin{cor}\label{Nil}
Every closed $\textrm{Nil}$ 3--manifold is diffeomorphic to a cusp
cross section of an arithmetic complex hyperbolic 2--manifold.
\end{cor}

\noindent Density for almost flat structures on compact
$\textrm{Nil}$ 3--manifolds follows from Theorem \ref{Complex} (a)
and Corollary \ref{Nil}. Our final result is the geometric
consequence of Corollary \ref{GenUnipotentSelbergLemma}.

\begin{cor}\label{MainGeometricCorollary}
Let $N$ be a closed infranil manifold and $X$ an arithmetic
orbifold. Then any proper $\pi_1$--injective immersion of $N$ into $X$ can be
lifted to be an embedding in a finite manifold cover of $X$.
\end{cor}

\noindent Note that Theorem \ref{SemisimpleAvoid} (see Theorem
\ref{SemisimpleArithmeticVersion}) permits one to find finite
manifold covers of $X$ such that any finite number of closed
geodesics fail to lift and $N$ can be lifted to be embedded.

\paragraph{Acknowledgements.}

\noindent I wish to thank Karel Dekimpe, Daniel Groves, Christopher Leininger, Gopal
Prasad, Alan Reid, Matthew Stover, and Henry Wilton for stimulating conversations
and gratefully acknowledge the California Institute of Technology
for their hospitality, as part of this work was done while visiting
that institution.

\section{Preliminaries}\label{Prelim}

\paragraph{Notation.} For each prime $p$, $\Z_p,\Q_p$ will denote the
$p$--adic integers and field, respectively. The full profinite closure of $\Z$
will be denoted by $\wh{\Z}$. Associated to these topological rings
are the topological groups $\GL(n,\Z_p), \GL(n,\Q_p)$, and
$\GL(n,\wh{\Z})$. Finally, $r_m\co \GL(n,\Z) \to \GL(n,\Z/m\Z)$ will denote the reduction
homomorphism given by reducing coefficients modulo $m$.

\paragraph{2.1.}\label{ProClose} Given a subgroup $\Gamma$ of $\GL(n,\Z)$, we
denote the closure of $\Gamma$ in $\GL(n,\Z_p)$ by $\Cl_p(\Gamma)$
and its closure in $\GL(n,\wh{\Z})$ by $\Cl(\Gamma)$. The following
is a restatement of \cite[Theorem 3.1]{McR1}.

\begin{thm}\label{Borel}
If $\Gamma<\GL(n,\Z)$ is virtually solvable, then $\Cl(\Gamma) \cap \GL(n,\Z) = \Gamma$.
\end{thm}

\paragraph{2.2.} Given an element $\ga$ in $\GL(n,\Z)$, there exists a unique
decomposition $\ga = \ga_s\ga_u$ called the \emph{Jordan
decomposition}. The elements $\ga_s,\ga_u \in \GL(n,\C)$ have the
following properties:

\begin{itemize}
\item[(1)]
$\ga_s$ is diagonalizable and $\ga_u-I_n$ is nilpotent.
\item[(2)]
$[\ga_s,\ga_u] \bdef \ga_s^{-1}\ga_u^{-1}\ga_s\ga_u=I_n$.
\end{itemize}

\noindent An element $\ga$ is called \emph{semisimple} if
$\ga_u=I_n$ and \emph{unipotent} if $\ga_s=I_n$.  It will be our
convention to consider the trivial element as unipotent. Whether or
not an element is semisimple or unipotent is a conjugacy invariant,
a fact gleamed from the formulae
\begin{equation}\label{ConjugacyFormulae}
 (\eta^{-1}\gamma\eta)_s = \eta^{-1}\ga_s\eta, \quad
(\eta^{-1}\gamma\eta)_u = \eta^{-1}\gamma_u\eta.
\end{equation}

\paragraph{2.3.} A subgroup $\Ga$ of $\GL(n,\C)$ is \emph{unipotent} if $\Ga$
is conjugate in $\GL(n,\C)$ into the group of upper triangular
matrices with ones along the diagonal. More generally, if $\Ga$ has
a finite index subgroup that is unipotent, we say that $\Ga$ is
\emph{virtually unipotent}. \smallskip\smallskip

\noindent Given a virtually unipotent subgroup $\Gamma$ of
$\GL(n,\Z)$, each element $\ga$ in $\Gamma$ possesses a Jordan
decomposition $\ga_s\ga_u$. As some power of $\ga$ is unipotent,
$\ga^m = \ga_u^m$ where $m$ is the order of $\ga_s$. In the event
that $\Gamma$ is torsion free, $\ga_u$ is necessarily nontrivial and hence
no element of $\Gamma$ can be semisimple. Note also that both $\ga_s,\ga_u$ reside in
$\GL(n,\Q)$.

\paragraph{2.4.} Associated to $\Gamma$ is the set of \emph{semisimple factors}
\[ \Semi(\Gamma) = \set{\ga_s~:~ \ga \in \Gamma} \su \GL(n;\C). \]
According to (\ref{ConjugacyFormulae}), the conjugate action of
$\Ga$ induces an action on the set $\Semi(\Ga)$. The finiteness of the quotient $\Semi(\Ga)/\Ga$ under this action will be critical.

\begin{lemma}\label{Finiteness}
If $\Gamma<\GL(n,\Z)$ is virtually unipotent, then
$\Semi(\Gamma)/\Ga$ is finite.
\end{lemma}

\noindent We postpone the proof of Lemma \ref{Finiteness} until
Section \ref{Tools}. For future reference, we fix a complete set of
representatives $^s\ga_1,\dots,~^s\ga_r \in \Semi(\Gamma)$ for the
quotient set $\Semi(\Gamma)/\Gamma$.

\section{Proof of Theorem \ref{SemisimpleAvoid}}\label{MainProof}

\noindent In this section we prove Theorem \ref{SemisimpleAvoid}. We begin by deducing Corollary \ref{UnipotentSelbergLemma}
from Theorem \ref{SemisimpleAvoid}.

\begin{pf}[Proof of Corollary \ref{UnipotentSelbergLemma}]
By Weil local rigidity, there are finitely many
conjugacy classes of torsion elements in $\GL(n,\Z)$ (see for instance
\cite{Rag}). Let $\eta_1,\dots,\eta_t$ be a complete set of
representatives for these conjugacy classes of torsion elements.
According to Theorem \ref{SemisimpleAvoid}, for each $\eta_j$, there
exists a finite index subgroup $\La_j$ of $\GL(n,\Z)$ such that
$\Ga<\La_j$ and $[\eta_j] \cap \La_j = \es$. The subgroup
\[ \La_0 = \ba_{j=1}^t \La_j, \]
is easily seen to suffice for verifying the corollary.
\end{pf}

\subsection{Proof of Theorem \ref{SemisimpleAvoid}}

\noindent In this subsection, we prove Theorem \ref{SemisimpleAvoid}. The proof is elementary (modulo Lemma \ref{Finiteness}), relying only Jordan form and passage to convergent subsequence (via compactness).

\paragraph{1. Some basic lemmas.} We begin by recording some elementary facts, the proofs of which
have been included for completeness.

\begin{lemma}\label{UnipotentLimits}
If $\ga$ is a limit of unipotent elements in $\GL(n,\Z_p)$, then
$\ga$ is unipotent.
\end{lemma}

\begin{pf}
Let $\set{\ga_j}$ be a sequence of unipotent elements in
$\GL(n,\Z_p)$ that converge to $\ga$. As there exists a uniform
bound on the multiplicative order of $\ga_j - I_n$, it follows that
$\ga$ is unipotent. Specifically, if $N$ is an integer such
that the multiplicative order of $\ga_j - I_n$ is bounded above by
$N$ for all $j$, it follows that for all $j>0$, $(\ga_j - I_n)^{N} =
0_n$. Thus
\begin{align*}
0_n &= \lim_j \pr{(\ga_j - I_n)^{N}} = \pr{\lim_j \pr{ \ga_j -
I_n}}^{N} = \pr{\pr{\lim_j \ga_j} - I_n}^{N} = (\ga -I_n)^{N}.
\end{align*}
\end{pf}

\begin{lemma}\label{ConjugacyClosures}
If $\eta \in \GL(n,\Z)$ is semisimple, then $\Cl_p([\eta])$ consists
of semisimple elements.
\end{lemma}

\begin{pf}
For $\la \in \Cl_p([\eta])$, there exists a convergent sequence
$\set{\eta'_j}$ in $[\eta]$ whose limit is $\la$. For each
$\eta_j'$, by definition there exists $\be_j \in \GL(n,\Z)$ such
that $\be_j^{-1}\eta_j'\be_j = \eta$. As $\set{\be_j}$ is a sequence
in the compact group $\GL(n,\Z_p)$, there exists a convergent
subsequence $\set{\be_\ell}$ of $\set{\be_j}$ with limit $\be \in
\GL(n,\Z_p)$. Note that by continuity of taking inverses, the sequence $\set{\be_\ell^{-1}}$ is also
convergent and has limit $\be^{-1}$. In total, this yields
\begin{align*}
\eta &= \lim_\ell \eta = \lim_\ell (\be_\ell^{-1} \eta_\ell' \be_\ell) = \pr{\lim_\ell \be_\ell^{-1}} \cdot \pr{\lim_\ell \eta_\ell'} \cdot \pr{\lim_\ell \be_\ell} = \be^{-1}\la\be.
\end{align*}
As $\eta$ is semisimple, $\la$ is as well.
\end{pf}

\begin{lemma}\label{CongruenceSeparation}
For subsets $R_1,R_2\su \GL(n,\Z)$, if $\Cl_p(R_1) \cap \Cl_p(R_2) =
\es$, then there exists a positive integer $K$ such that $r_{p^K}(R_1) \cap
r_{p^K}(R_2) = \es$.
\end{lemma}

\begin{pf}
Note that as the closed sets $\Cl_p(R_1)$ and $\Cl_p(R_2)$ are
disjoint, the topological normality of $\GL(n,\Z_p)$ implies that we
can find an open subsets $O_j$ of $\GL(n,\Z_p)$ that contain
$\Cl_p(R_j)$ and are disjoint from $\Cl_p(R_k)$ where $j\ne k$. The
subsets $\Cl_p(r_{p^\ell}^{-1}(r_{p^\ell}(R_j)))$ are open (and closed) in
$\GL(n,\Z_p)$, contain $\Cl_{p}(R_j)$, and have the feature that
\[ \ba_{\ell=1}^\iny \Cl_{p}(r_{p^\ell}^{-1}(r_{p^\ell}(R_j))) = \Cl_{p}(R_j). \]
Therefore, for some large integer $K$, it must be that
\[ \Cl_p(r_{p^K}^{-1}(r_{p^K}(R_1))) \cap \Cl_p(r_{p^K}^{-1}(r_{p^K}(R_2))) = \es. \]
Thus, we must have the less restrictive, desired conclusion
\[ r_{p^K}(R_1) \cap r_{p^K}(R_2) = \es. \]
\end{pf}

\paragraph{2. Limit point criterion.} For the statement of the following proposition, recall by Lemma \ref{Finiteness} that there exists a finite set
$\set{~^s\ga_1,\dots,~^s\ga_r}$ of semisimple factors up to $\Gamma$--conjugation.

\begin{prop}\label{LimitCriterion}
If $\eta \in \Cl_p(\Gamma)$ is semisimple, then there exists
$1\leq k_\eta\leq r$ such that $^s\ga_{k_\eta} \in \Cl_p(\Gamma)$.
\end{prop}

\begin{pf} By definition, there exists a convergent sequence $\set{\ga_j}$ in
$\Gamma$ with limit $\eta$. Consider the pair of sequences $^s\ga_j = (\ga_j)_s$, $^u\ga_j = (\ga_j)_u$. We will first prove the proposition under the assumption that $^s\ga_j = ~^s\ga_{k_\eta}$ for all $j$ and some fixed $k_\eta$. We will see below that the general situation can be reduced to this. Under the assumption that $^s\ga_j =~^s\ga_{k_\eta}$, the
associated unipotent factor sequence $\set{~^u\ga_j}$ of $\set{\ga_j}$ must also
converge since $^u\ga_j = \ga_j(~^s\ga_{k_\eta})^{-1}$. Suggestively
setting $\eta_u$ to be the limit of the sequence $\set{~^u\ga_j}$,
we assert that $^s\ga_{k_\eta}\eta_u$ is the Jordan decomposition
for $\eta$. That $\eta_u$ is unipotent follows from Lemma
\ref{UnipotentLimits} (we already know that $^s\ga_{k_\eta}$ is
semisimple). To see that $^s\ga_{k_\eta}\eta_u = \eta$, notice that
$^s\ga_{k_\eta}\cdot ~^u\ga_j = \ga_j$. Therefore,
\begin{align*}
\eta &= \lim_j \ga_j = \lim_j \pr{~^s\ga_{k_\eta} \cdot ~^u\ga_j} = ~^s\ga_{k_\eta} \cdot \pr{\lim_j~^u\ga_j} = ~^s\ga_{k_\eta}\cdot \eta_u.
\end{align*}
Finally to see that $[~^s\ga_{k_\eta},\eta_u]=I_n$, note that
\begin{align*}
I_n &= \lim_j [~^s\ga_{k_\eta},~^u\ga_j] =
\lim_j \pr{(~^s\ga_{k_\eta})^{-1}\cdot (~^u\ga_j)^{-1} \cdot ~^s\ga_{k_\eta}\cdot ~^u\ga_j} \\
&= (~^s\ga_{k_\eta})^{-1}\cdot \pr{\lim_j (~^u\ga_j)^{-1}} \cdot ~^s\ga_{k_\eta} \cdot \pr{\lim_j ~^u\ga_j} \\
&= (~^s\ga_{k_\eta})^{-1}\cdot \eta_u^{-1} \cdot ~^s\ga_{k_\eta}
\cdot \eta_u = [~^s\ga_{k_\eta},\eta_u]
\end{align*}
as needed. This shows that $^s\ga_{k_\eta}\eta_u$ is the Jordan
decomposition for $\eta$. As $\eta$ is semisimple, it must be that
$\eta_u=I_n$ and hence $\eta = ~^s\ga_{k_\eta}$ for some $k_\eta$.\smallskip\smallskip

\noindent It could be the case that the semisimple factor sequence $^s\ga_j$ for $\ga_j$ is not constant. Using Lemma \ref{Finiteness}, we will reduce this case to the previous one. To begin, by Lemma \ref{Finiteness}, there exists a sequence $\set{\al_j}$ in
$\Gamma$ such that
\[ (\al_j^{-1}\ga_j\al_j)_s = \al_j^{-1}~^s\ga_j \al_j =
~^s\ga_{k_j}, \quad k_j \in \set{1,\dots,r}. \]
In particular, some $k_\eta$ must occur infinitely often and so we can
pass to a subsequence $\ga_i$ such that
\[ (\al_i^{-1}\ga_i\al_i)_s = ~^s\ga_{k_\eta} \]
for some fixed $1 \leq k_\eta \leq r$. As $\set{\al_i}$ is a sequence in the
compact group $\Cl_p(\Gamma)$, $\set{\al_i}$ has a convergent
subsequence $\set{\al_\ell}$ with limit $\al \in \Cl_p(\Ga)$.
Again by continuity of taking inverses, $\set{\al_\ell^{-1}}$ is convergent with limit $\al^{-1} \in \Cl_p(\Gamma)$. In total, we see now that
\begin{align*}
\lim_\ell \al_\ell^{-1} \ga_\ell \al_\ell
&= \pr{\lim_\ell \al_\ell^{-1}} \cdot \pr{\lim_\ell \ga_\ell} \cdot \pr{\lim_\ell \al_\ell} = \al^{-1} \eta \al.
\end{align*}
As $\al,\al^{-1},\eta \in \Cl_p(\Gamma)$, so is $\al^{-1}\eta\al$. In
addition, since $\eta$ is semisimple, so is its conjugate $
\al^{-1}\eta\al$. By taking $\al^{-1}\eta \al$ instead of $\eta$, we
can assume that $\eta$ is the limit of a sequence $\set{\ga_j}$ in
$\Gamma$ whose semisimple factors are constant.
\end{pf}

\paragraph{3. Avoiding a semisimple factor.} As before, the
elements $^s\ga_1,\dots,~^s\ga_r$ are a complete list of semisimple
factors up to $\Ga$--conjugation given by Lemma \ref{Finiteness}.

\begin{lemma}\label{UnipotentClosure}
For each $k=1,\dots,r$, there exists a prime $p_k$ such that
$^s\ga_k \notin \Cl_{p_k}(\Gamma)$.
\end{lemma}

\begin{pf}
If $^s\ga_k \notin \GL(n,\Z)$, then there exists a matrix coefficient $\nu$
of $^s\ga_k$ such that $\nu \notin \Z$. Taking $p_k$ to be a prime
occurring in the denominator of $\nu$, it follows that $\nu \notin \Z_{p_k}$.
As any limit of elements in $\Gamma$ is in $\GL(n,\Z_p)$, $^s\ga_k
\notin \Cl_{p_k}(\Gamma)$. We now consider the alternative when $^s\ga_k \in \GL(n,\Z)$. According to
Theorem \ref{Borel}, if $^s\ga_k \in \GL(n,\Z) \cap \Cl(\Gamma)$,
then $^s\ga_k \in \Gamma$. However, $\Gamma$ is torsion free and
$^s\ga_k$ is finite order, and thus could not possibly reside in $\Ga$. Therefore, there must exist a prime $p_k$ such that $^s\ga_k \notin \Cl_{p_k}(\Gamma)$, as desired.
\end{pf}

\paragraph{4. Proof of Theorem \ref{SemisimpleAvoid}.} Let $[\eta]$ be a $\GL(n,\Z)$--conjugacy class for a semisimple element
$\eta$ in $\GL(n,\Z)$. Using the primes in Lemma
\ref{UnipotentClosure} and setting
\[ N = \prod_{i=1}^r p_i, \quad \Cl_N(\Gamma) = \prod_{i=1}^r \Cl_{p_i}(\Gamma), \]
it follows that $^s\ga_k \notin \Cl_N(\Gamma)$ for all $k=1,\dots r$. In particular,
$\Cl_N(\Gamma)$ contains no semisimple elements by
Proposition \ref{LimitCriterion}. By Lemma \ref{ConjugacyClosures},
$\Cl_N([\eta])$ consists entirely of semisimple
elements. These two facts imply that $\Cl_N(\Ga) \cap
\Cl_N([\eta]) = \es$. By Lemma \ref{CongruenceSeparation}, there exists a positive integer $K$ such that
$r_{N^K}(\Ga) \cap r_{N^K}([\eta]) = \es$. The proof is completed by taking the finite
index subgroup $r_{N^K}^{-1}(r_{N^K}(\Ga))$ for $\Gamma_0$.\qed\smallskip\smallskip\smallskip

\noindent Theorem \ref{SemisimpleAvoid} is the strongest possible
result. In Section \ref{Caution}, we give an example, due to Stebe \cite{Stebe}, of an infinite cyclic
subgroup of $\GL(n,\Z)$ with semisimple generator for which Theorem \ref{SemisimpleAvoid} is
false. In particular, the virtual unipotency assumption cannot not
be dropped.

\subsection{Torsion in profinite groups}

For a torsion free, residually finite $G$, there is no reason
to expect the profinite closure $\widehat{G}$ of $G$ to be torsion
free. Indeed, torsion free, finite index subgroups of $\GL(n,\Z)$ with $n>2$ provide linear
examples (see \cite{Lub}). Even for nilpotent groups $G$, it need
not be the case that $\widehat{G}$ is torsion free (see \cite{KW}).
However, for the class of $\Gamma$ consider here, it follows from
\cite{KW} that $\wh{\Gamma}$ is torsion free. In addition, it
follows from \cite{McR1} that $\Cl(\Gamma) = \wh{\Gamma}$. Using
this with Lemma \ref{ConjugacyClosures} provides a different proof
of Corollary \ref{UnipotentSelbergLemma}. Our proof of
Theorem \ref{SemisimpleAvoid} provides an elementary proof that $\wh{\Gamma}$ is
torsion free for virtually unipotent subgroups of $\GL(n,\Z)$.

\section{Proof of Lemma \ref{Finiteness}}\label{Tools}

\noindent In this section, we prove Lemma \ref{Finiteness}. We refer
the reader to \cite{Dek} for the material used below on nilpotent
Lie groups, Lie algebras, and almost crystallographic groups.

\paragraph{Preliminaries.} For a virtually unipotent subgroup $\Ga$ of $\GL(n,\Z)$, there exists a short exact sequence
\[ 1 \lra \Ga_u \lra \Ga \lra \te \lra 1 \]
where $\Ga_u$ is the Fitting subgroup of $\Ga$ and $\te$ is a finite group (the holonomy group of $\Ga$). The associated
holonomy representation $\vp\co \te \to \Out(\Ga_u)$ together with a
2--cocycle $f \in H_\vp^2(\Ga_u,\te)$ determine $\Ga$. We will prove
Lemma \ref{Finiteness} by induction of the step size of $\Ga_u$. The
base case when $\Ga_u$ is abelian is nothing more than the case when
$\Ga$ is a crystallographic group. Before addressing the base case,
we simplify our situation.\smallskip\smallskip

\noindent Set $\B{N}$ to be the Mal'cev completion \cite[p. 9]{Dek}
of $\Ga_u$ and $\Fr{n}$ to be the Lie algebra of $\B{N}$. By
construction, $\Ga_u$ admits an injection into $\B{N}$.
The group $\B{N}$ is a connected, simply connected, nilpotent Lie group and so the exponential map (see
\cite[p. 7--8]{Dek}) $\exp\co \Fr{n} \to \B{N}$ has a smooth inverse
$\log\co \B{N} \to \Fr{n}$. By Mal'cev rigidity \cite[Theorem 1.2.3]{Dek}, the
holonomy representation $\vp$ has a unique extension $\ol{\vp}\co
\te \to \Out(\B{N})$ and this extension lifts to a homomorphism into
$\Aut(\B{N})$ (see \cite[Lemma 3.1.2]{Dek}). This provides us with
an injection $\psi\co \Ga \to \B{N} \rtimes_{\ol{\vp}} \te$ where,
in an abuse of notation, $\ol{\vp}$ denotes some lift of $\ol{\vp}$
to $\Aut(\B{N})$. This allows us to write each element $\ga \in \Ga$
as $(n_\ga,\te_\ga)$ with $n_\ga \in \B{N}$ and $\te_\ga \in \te$.
We also have a Jordan decomposition of $\ga$ given by $\ga =
(n_s,\te_\ga)\cdot (n_u,1)$ where $n_s,n_\ga \in \B{N}$ and
$\te_\ga(n_u)=n_u$. The set of semisimple factors under this
decomposition is given by
\[ \Semi_{\B{N}}(\Ga) = \set{ (n_s,\te_\ga)~:~\ga \in \Ga} \su \B{N} \rtimes_{\ol{\vp}} \te. \]
and we can reduce the finiteness of $\Semi(\Ga)/\Ga$ to the
finiteness of $\Semi_{\B{N}}(\Ga)/\Ga$. That this can be done is
seen by the following argument. By Mal'cev rigidity, the inclusion
of $\Ga$ into $\GL(n,\Z)$ induces a smooth injection $\rho\co
\B{N}\rtimes_{\ol{\vp}} \te \to \GL(n,\R)$. By the uniqueness of the
Jordan decomposition (see \cite[I.4]{Borel1}), we have that
$\rho((n_s,\te_\ga)) = \ga_s$, $\rho((n_u,1)) = \ga_u$.
Consequently, it suffices to show the finiteness of
$\Semi_{\B{N}}(\Ga)/\Ga$. We are now ready to prove Lemma \ref{Finiteness}.

\paragraph{Proof of Lemma \ref{Finiteness}.}
Our proof will be done by inducting on the step size of $\Ga_u$.\smallskip\smallskip

\noindent \textbf{Base case.} In this case $\Ga_u \cong
\Z^m$ for some $m$ and hence $\B{N}=\R^m$. By the Bieberbach
theorems (see \cite{Ch}), we write elements as $(t,S)$ where $t \in
\Z^m$ and $S \in \GL(m,\Z)$. As there are only finitely many $S$
(these are the elements of $\te$), it suffices to prove that there
are only finitely many semisimple factors $(t_s,S)$ up to
$\Ga$--conjugation for each individual $S$. The action of $S$ on
$\Q^m$ decomposes into two subspaces $\Q^m = W_S \op W_{triv,S}$
where $W_{triv,S}$ is the maximal subspace of $\Q^m$ on which $S$
acts trivially. It is a simple matter to see that the Jordan
decomposition of an element $(t,S)$ is of the form $(t_s,S)(t_u,I_m)$
where $t_s \in W_S$ and $t_u \in W_{triv,S}$. Conjugating by
$(t,I_m)$ produces $(t_s + (S-I_m)t,S)$. As we are only concerned
with those vectors in $W_S$, we may assume $t \in W_S$. The possible
vectors $t$ form a finite index $\Z$--submodule of $W_S(\Z)$ whose
image under $S-I_m$ is still a finite index $\Z$--submodule of
$W_S(\Z)$ since $S-I_m$ is invertible on $W_S$. As the set of
possible vectors $t_s$ is contained in $W_S(\Z)$, up to
$\Ga_u$--conjugacy, the possible vectors are identified with a
subset of the quotient $W_S(\Z)/(S-I_m)(L)$, where $L$ is the
$\Z$--submodule of vectors in $W_S(\Z)$ which arise as translation
vectors for an element of $\Ga_u$.  As this quotient is finite, we
conclude $\Semi_\B{N}(\Ga)/\Ga_u$ is finite and thus $\Semi_\B{N}(\Ga)/\Ga$ is finite.\smallskip\smallskip

\noindent \textbf{General case.} For the general case, let $\Ga_u^k$
denote the $k$th term in the lower central series for $\Ga_u$ and
assume that $\Ga_u$ has step size $j>1$ (i.e., $\Ga_u^j = \set{1}$).
Associated to each $\Ga_u^k$ is its Mal'cev completion $\B{N}_k$ and
Lie algebra $\Fr{n}_k$. The conjugate action of $\Ga$ on $\B{N} \rtimes_{\ol{\vp}}
\te$ induces an $\Ad(\Ga)$--action on $\Fr{n} \rtimes_{\ol{\vp}}
\te$. The semisimple factor set
$\Semi_{\B{N}}(\Ga)$ produces a corresponding set
\[ \Semi_{\Fr{n}}(\Ga) = \set{ (\eta_s,\te_\ga)~:~ \ga \in
\Ga,~\eta_s = \log(n_s)} \su \Fr{n} \rtimes_{\ol{\vp}} \te. \] The finiteness of
$\Semi_{\B{N}}(\Ga)/\Ga$ is equivalent to the finiteness of
$\Semi_{\Fr{n}}(\Ga)/\Ad(\Ga)$. Consequently, it
suffices to show the latter. In addition, it suffices to show the
finiteness of $\Semi_{\Fr{n}}(\Ga)/\Ad(\Ga_u)$
as $\abs{\Semi_{\Fr{n}}(\Ga)/\Ad(\Ga_u)}$ is at least as big as
$\abs{\Semi_{\Fr{n}}(\Ga)/\Ad(\Ga)}$. We will now establish the finiteness of the latter set as follows. The Lie algebra $\Fr{n}$ of $\B{N}$ is a graded vector space of the form
\[ \Fr{n} = \bop_{i=0}^{j-1} \Fr{n}_i/\Fr{n}_{i+1} = \bop_{i=0}^{j-1} \textrm{Gr}_i(\Fr{n}) \]
where $\Fr{n}_0 = \Fr{n}$ and $\Fr{n}_j = \set{0}$. In particular,
each element $\eta_s$ has the form
\begin{equation}\label{GradedForm}
\eta_s = (\eta_0,\dots,\eta_{j-1}), \quad \eta_i \in
\textrm{Gr}_i(\Fr{n}).
\end{equation}
Notice that we have a pair of almost crystallographic groups $\Ga'$
and $\Ga''$ given by
\[ 1 \lra \Ga_u^{1} \lra \Ga' \lra \te \lra 1 \]
and
\[ 1 \lra \Ga_u/\Ga_u^{1} \lra \Ga'' \lra \te \lra 1. \]
This pair of groups inject into the groups $\B{N}_{1}
\rtimes_{\ol{\vp}} \te$ and $(\B{N}/\B{N}_{1}) \rtimes_{\ol{\vp}}
\te$, respectively. For $\Ga'$, we have an induced
$\Ad(\Ga_u^{1})$--action on $\Fr{n}_{1} \rtimes_{\ol{\vp}} \te$
where the latter is nothing more than
\[ \Fr{n}_{1} = \bop_{i=1}^{j-1} \Fr{n}_{i}/\Fr{n}_{i+1} = \bop_{i=1}^{j-1} \textrm{Gr}(\Fr{n}). \]
Likewise, we have an $\Ad(\Ga_u/\Ga_u^{1})$--action on
$(\Fr{n}/\Fr{n}_{1}) \rtimes_{\ol{\vp}} \te$. According to our
induction hypothesis, there only finitely many possibilities for
$\eta_1,\dots,\eta_{j-1}$ in (\ref{GradedForm}) up to the
$\Ad(\Ga_u^{1})$--action. Similarly, by our induction hypothesis,
there are only finitely many possibilities for $\eta_0$ in
(\ref{GradedForm}) up to the $\Ad(\Ga_u/\Ga_u^{1})$--action. Thus,
there are only finitely many possibilities for
$\eta_0,\dots,\eta_{j-1}$ in (\ref{GradedForm}) up to the
$\Ad(\Ga_u)$--action. In particular, up to the $\Ad(\Ga_u)$--action,
there are only finitely many possibilities for $\eta_s$ in
$(\eta_s,\te_\ga) \in \Semi_{\Fr{n}}(\Ga)$. As the possibilities for
$\te_\ga$ range over the finite group $\te$, this implies the
finiteness of $\Semi_{\Fr{n}}(\Ga)/\Ga$.\qed

\section{Theorem \ref{SemisimpleAvoid} for arithmetic lattices}\label{GenSect}

\noindent The proof of Theorem \ref{SemisimpleAvoid} and its Corollary
\ref{UnipotentSelbergLemma} work for subgroups $\tri$ of $\GL(n,\Q)$
commensurable with $\GL(n,\Z)$. Briefly we describe this
here. We begin with the following lemma whose validity can be deduced from the proof that $\Ga$ injects into $\B{N} \rtimes_{\ol{\vp}} \te$.

\begin{lemma}\label{TorsionLattice}
There exists a lattice $\Ga_0 <\B{N} \rtimes_{\ol{\vp}} \te$ such that $\Ga_0$
contains each $^s\ga_k$ and $\Ga$.
\end{lemma}

\noindent With Lemma \ref{TorsionLattice}, we can generalize Theorem
\ref{SemisimpleAvoid}. To this end, let $\tri$ be a
subgroup of $\GL(n,\Q)$ commensurable with $\GL(n,\Z)$ and assume
that $\tri$ contains a torsion free, virtually unipotent subgroup
$\Ga$. Using the same approach as in the proof of
Theorem \ref{SemisimpleAvoid}, note that Proposition
\ref{LimitCriterion} is validated as before (note that passing to convergent subsequences
is done now inside the compact set $\textrm{Cl}_p(\tri)$). For Lemma
\ref{UnipotentClosure}, we must modify our argument. It could be the case that $\tri$ does not contain the
elements $^s\ga_k$ coming from Lemma \ref{Finiteness}. However, by
Lemma \ref{TorsionLattice} and \cite[Corollary 10.14]{Rag}, there
exists a group $\tri_0$ commensurable with $\tri$ that contains
$\Ga_0$. By Theorem \ref{Borel} (this holds for groups commensurable with $\GL(n,\Z)$), $\Cl(\Ga) \cap \tri_0 = \Ga$. In
particular, for each $^s\ga_k$, there must exist a prime $p_k$ such
that $^s\ga_k \notin \Cl_{p_k}(\Ga)$. This shows that Theorem
\ref{SemisimpleAvoid} can be extended to groups $\tri$ in
$\GL(n,\Q)$ commensurable with $\GL(n,\Z)$. For a general arithmetic lattice $\La$, there exists an injective homomorphism $\psi\co \La \to \GL(n,\Q)$ such that $\psi(\La)$ is contained in a subgroup $\tri$ in $\GL(n,\Q)$ that is commensurable with
$\GL(n,\Z)$. Using the above argument, for any semisimple element $\eta \in \La$, we can
find a finite index subgroup $\tri_0$ of $\tri$ such that
$\Ga<\tri_0$ and $[\eta]_{\tri} \cap \tri_0 = \es$ where
$[\eta]_\tri$ is the $\tri$--conjugacy class of $\eta$. Certainly
$[\eta]_\La \su [\eta]_\tri$ and thus $[\eta]_\La \cap \tri_0=\es$.
Intersecting $\tri_0$ with $\La$ produces a finite index subgroup
$\La_0$ of $\La$ such that $\Ga<\La_0$ and $[\eta]_\La \cap \La_0 =
\es$. In total, we obtain the following theorem and corollary.

\begin{thm}\label{SemisimpleArithmeticVersion}
Let $\La$ be an arithmetic lattice, $\Ga<\La$ a torsion free,
virtually unipotent subgroup, and $\eta \in \La$ a semisimple
element. Then there exists a finite index subgroup $\La_0$ of $\La$
such that $\Ga<\La_0$ and $[\eta]_\La \cap \La_0 = \es$.
\end{thm}

\begin{cor}\label{GenUnipotentSelbergLemma} Let $\La$ be an
arithmetic lattice and $\Ga<\La$ a torsion free, virtually unipotent
subgroup. Then there exists a torsion free finite index subgroup
$\La_0$ of $\La$ such that $\Ga<\La_0$.
\end{cor}

\noindent The arithmetic assumption is only used in the proof of Lemma \ref{UnipotentClosure}.
Thus, we have the following corollary.

\begin{cor}
Let $\La<\GL(n,\C)$ be a finitely generated group and $\Ga<\La$ a
torsion free, virtually unipotent subgroup. Given an infinite
order semisimple element $\eta \in \La$, there exists a finite index
subgroup $\La_0<\La$ such that $\Ga<\La_0$ and $[\eta]_\La \cap \La_0 =
\es$.
\end{cor}

\noindent This corollary follows from the fact
that any semisimple $\eta$ in $\Cl_p(\Gamma)$ is conjugate to a
torsion element and thus itself is torsion. Indeed, there is nothing
special about taking the conjugacy class of an infinite order semisimple element.
The following is a consequence of the same logic.

\begin{cor}
Let $\La<\GL(n,\C)$ be a finitely generated group, $\Ga<\La$ a torsion free, virtually
unipotent subgroup, and $C$ be an infinite cyclic subgroup generated by a semisimple element.
Then $\Cl(\Gamma) \cap \Cl(C)=\set{1}$.
\end{cor}

\noindent Wilton--Zalesskii have also obtained this result in the case when $\La$ is a Kleinian group and $\Ga$ is a parabolic subgroup.\smallskip\smallskip

\noindent \textbf{Remark.} We mention in passing that one can prove Corollary \ref{GenUnipotentSelbergLemma} as before using the fact that $\wh{\Ga}$ is torsion free, $\Cl(\Ga) = \wh{\Ga}$, and Lemma \ref{ConjugacyClosures}.

\section{Geometric applications of Corollary
\ref{GenUnipotentSelbergLemma}}\label{GeometricAppsSect}

\noindent In this section, we derive the main geometric corollaries
of Corollary \ref{GenUnipotentSelbergLemma} mentioned in the
introduction. For brevity, we refer the reader to
\cite{LR,McR1,McR2} for some of the details.

\subsection{Flat manifolds}

\paragraph{1. Proof of Corollary \ref{Arith}.}

\noindent For a fixed flat $n$--manifold $N$, Long and Reid \cite{LR}
constructed an arithmetic hyperbolic $(n+1)$--orbifold $M$ such that
$N$ is diffeomorphic to a cusp cross section of $M$. In particular,
$\pi_1(N)$ is a torsion free, virtually unipotent subgroup of $\pi_1^{orb}(M)$. By
Corollary \ref{GenUnipotentSelbergLemma}, there exists a finite
index, torsion free subgroup $\La_0<\pi_1^{orb}(M)$ such that
$\pi_1(N)<\La_0$. Passing to the cover $M_0 \to M$ corresponding to
$\La_0$, yields an arithmetic hyperbolic $(n+1)$--manifold $M_0$
such that $N$ is diffeomorphic to a cusp cross section of $M_0$.\qed

\paragraph{2. Nimershiem's conjecture.}

\noindent Reviewing the proof of Corollary \ref{Arith}, notice that
passage from $M$ to $M_0$ does not change the flat similarity class
on the cusp cross section diffeomorphic to $N$. In particular, we
obtain the following orbifold-to-manifold promotion.

\begin{cor}\label{OrbifoldToManifold}
The space of flat similarity classes on a flat $n$--manifold that
arise in cusp cross sections of arithmetic hyperbolic
$(n+1)$--orbifolds is precisely the same as those that arise in
arithmetic hyperbolic $(n+1)$--manifolds.
\end{cor}

\noindent We established \cite[Proposition 3.2]{McR2} the density of
those similarity classes that arise in cusp cross sections of
arithmetic orbifolds. This with Corollary \ref{OrbifoldToManifold}
proves Theorem \ref{Nimershiem}.

\paragraph{3. Classification of arithmetic cusp shapes.}

\noindent One of the main motivations for the geometric results of
this article come from Gromov \cite{Gromov} whose worked inspired
the conjectures of Farrell--Zdravkovska \cite{FarrellZdravkovska83}
and Nimershiem \cite{Nim1}. In the former, it was conjectured that
every flat $n$--manifold was diffeomorphic to the cusp cross section
of a one cusped hyperbolic $(n+1)$--manifold. However, Long and Reid
\cite{LongReid00} found examples of flat $3$--manifolds that can
never be diffeomorphic to a cusp cross section of a one cusped
hyperbolic 4--manifold. Corollary \ref{Arith} shows the conjectural
picture proposed by Farrell--Zdravkovska is not too far off (in some sense).
\smallskip\smallskip

\noindent Corollary \ref{OrbifoldToManifold} and \cite[Theorem
3.7]{McR2} show the set of flat similarity classes appearing as cusp
shapes in arithmetic hyperbolic manifolds is the image of the
rational points of an algebraic set under a projection map. In
total, this classifies cusp shapes of arithmetic hyperbolic
$(n+1)$--manifolds.

\subsection{Infranil manifolds}

\noindent The generalizations for complex and quaternionic
hyperbolic spaces, namely Theorem \ref{Complex}, follows from an
identical argument using Corollary \ref{GenUnipotentSelbergLemma}
and \cite[Theorem 3.5]{McR2}. The density of these structures in the
case $N$ is a $\Nil$ 3--manifold follows from Corollary
\ref{GenUnipotentSelbergLemma} and \cite[Corollary 3.6]{McR2}.
Finally, Corollary \ref{MainGeometricCorollary} follows from
Corollary \ref{GenUnipotentSelbergLemma} and \cite[Theorem
3.12]{McR1}. Just as there are flat $n$--manifolds that cannot arise as
the cusp cross section of a single cusped hyperbolic $(n+1)$--manifold,
there exist $\Nil$ 3--manifolds that cannot arise as the cusp cross
section of a one cusped complex hyperbolic 2--manifold (see \cite{Kamishima}). Corollary
\ref{MainGeometricCorollary} again shows the failure is not total.

\section{Final remarks}

\paragraph{1. Generalizing Theorem \ref{SemisimpleArithmeticVersion}.}\label{Caution}

\noindent For a virtually unipotent, torsion free subgroup $\Ga$,
there is essentially no difference in separating $\Ga$ from a
semisimple class or a torsion class. Even for an infinite cyclic
group $\innp{A}$ generated by a semisimple element $A$, it could
very well be the case that one cannot separate $\innp{A}$ from a
fixed semisimple conjugacy class $[B]$. Indeed, the failure of
conjugacy separability in $\SL(n,\Z)$, $n>2$ provides examples (see
\cite{Stebe}). However, the elements $A$ and $B$ are conjugate in
$\SL(n,\C)$ ($A$ and $B$ are conjugate in $\SL(n,\Z_p)$) and
thus it is possible to separate $\innp{A}$ from a fixed
torsion class. Indeed, using Theorem \ref{Borel}, Lemma \ref{ConjugacyClosures},
and the fact that $\wh{\innp{A}}$ is torsion free, one can find a torsion free
finite index subgroup of $\SL(n,\Z)$ that contains $\innp{A}$. In fact,
when $A$ is semisimple, this does not require an arithmetic assumption either.

\paragraph{2. Higher rank cusp cross sections.}

\noindent For cusp cross sections of higher rank locally symmetric
spaces, the fundamental group of a cusp cross section is virtually
solvable but typically not virtually unipotent. For instance, cusp
cross sections of Hilbert modular surfaces are $\textrm{Sol}$
3--orbifolds (see \cite{McR3} for more on this). Though Theorem \ref{SemisimpleAvoid} might not hold for
these groups, Corollary \ref{UnipotentSelbergLemma} extends. Indeed,
the profinite completion of such torsion free groups are known to be torsion free
by \cite{KW} and the profinite completion is isomorphic to
$\Cl(\Gamma)$ by \cite{McR1}. This with Lemma
\ref{ConjugacyClosures} implies Corollary
\ref{UnipotentSelbergLemma} for these groups.

\paragraph{3. Totally geodesic, immersed surfaces.}

\noindent In general, it seems difficult to resolve torsion in Question 2 from the introduction even
when $M$ is a hyperbolic 3--orbifold and $N$ is a totally geodesic
surface. However, there are some special cases when this can be
done. Indeed, when $\pi_1(M)$ is subgroup separable, since
$\wh{\pi_1(N)}$ is torsion free and the closure of $\pi_1(N)$ in
$\wh{\pi_1(M)}$ is isomorphic to $\wh{\pi_1(N)}$, one can extend
Corollary \ref{UnipotentSelbergLemma}. One class of $M$ that satisfy
this condition are noncompact arithmetic hyperbolic 3--orbifolds
(see \cite{ALR}) which are endowed with many totally geodesic,
immersed surfaces (see \cite{MR}).



\noindent Department of Mathematics \\
University of Chicago \\
Chicago, IL 60637 \\
email: {\tt dmcreyn@math.uchicago.edu}


\end{document}